
\magnification=\magstep1
\vsize=23.5true cm
\hsize=16true cm
\parskip=2pt
\topskip=1true cm
\headline={\tenrm\hfil\hfil}
\raggedbottom
\abovedisplayskip=3mm
\belowdisplayskip=3mm
\abovedisplayshortskip=0mm
\belowdisplayshortskip=2mm
\normalbaselineskip=15pt
\normalbaselines


\def\Z   {{\bf Z}}

\def\R   {{\bf R}}
\def\C   {{\bf C}}
\def\P   {{\bf P}}

\def\F   {{\bf F}}


\def\rank  {\mathop{\rm rank}}

\def\dim   {{\rm dim}}


\def\HZ      #1#2{{\rm H}^{#1}(#2,\Z)}
\def\HZhom   #1#2{H_{#1}(#2,\Z)}

\def\HR      #1#2{H^{#1}(#2,\R)}

\def\h       #1#2#3{h^{#1}(#2,#3)}

\def\hOE     #1{h^{#1}({\cal O} _E)}


\def\Kbar      {\overline{\cal K}}

\def\Kahler    {K\"ahler\ }
\def\CY        {Calabi--Yau\ }

\def\Xbar      {\bar X} 
\def\phimapX   {\phi : X \to \Xbar }

\def\tE        {{\tilde E}}
\def\tC        {{\tilde C}}
\def\3f        {{\CY threefold\ }}

\def\phrc      {{pseudo-holomorphic rational curve }}

\def\curves    {{pseudo-holomorphic rational curves }}
\def\rdp       {{rational double point }}

\mathchardef\checkchar"7014


\def\calX      {{\cal X}}
\def\calO      {{\cal O}}

\def\calM      {{\cal M}}
\def\calU      {{\cal U}}

\def\calJ      {{\cal J}}

\def\calK      {{\cal K}}


\def\BW        {1}
\def\CFKM      {2}
\def\Dem       {3}
\def\Fr        {4}
\def\GrCan     {5}
\def\GrDef     {6}
\def\HW        {7}
\def\Kont      {8}
\def\Loo       {9}
\def\McD       {10}
\def\McDS      {11}
\def\RCan      {12}
\def\RMin      {13}
\def\ReDP      {14}
\def\Ru        {15}
\def\Ruan      {16}
\def\RT        {17}
\def\Sik       {18}
\def\Wall      {19}
\def\WPic      {20}
\def\WKc       {21}
\def\WEll      {22}

\bigskip
\centerline{\bf SYMPLECTIC DEFORMATIONS OF CALABI--YAU THREEFOLDS}
\bigskip  \centerline{\bf P.M.H. WILSON} \bigskip\bigskip 
\bigskip \centerline{\bf Introduction} \bigskip \rm

Much recent work on \CY threefolds has concentrated on the analogy between
small deformations of the complex structure on a \3f $X$, which by the
unobstructedness theorem of Bogomolov, Tian and Todorov are parametrized by a
neighbourhood of the origin in $H^1 (T_X ) = H^{2,1} (X, \C )$, and small
deformations of the symplectic structure, which by Moser's Theorem are
parametrized by a  neighbourhood of the origin in $H^2 (X, \R ) = H^{1,1} (X, \R
)$.  This analogy is described in [\Kont ] in the context of Mirror Symmetry.

We recall that a \3f has various differential invariants, in particular the
cubic form on $\HZ 2 X$ given by cup-product, the linear form on $\HZ 2 X$
given by cup-product with the first Pontryagin class $p_1$ (which in the \CY
case is just the second chern class $c_2$), and the middle cohomology $\HZ 3 X$.
Indeed, if $X$ is simply connected and has no torsion in $\HZ 3 X$, these
invariants specify the differential class of the manifold precisely [\Wall ].

In this paper, we shall denote a \3f by a pair $(X,J)$, where $X$ is the
differential manifold and $J$ is the (integrable) complex structure.  A
symplectic form $\omega $ on $X$ which is compatible with $J$ is just a \Kahler
form on the complex manifold.  Moreover, the almost complex structures
which are compatible with a given symplectic form on $X$ form a contractible
set, as do the almost complex structures tamed by the symplectic form.  It
follows therefore that if we have a family of symplectic forms $\omega _t$ on
$X$, we can find a corresponding family $J_t$ of almost complex structures, with
$J_t$ tamed by $\omega _t$.  We are led therefore to considering \it almost
complex \rm deformations of $J$ -- for very accessible surveys of the concepts
from Symplectic Geometry we need, the reader is referred to [\McD , \McDS ].  It
will follow from our results that if two \CY threefolds (with given \Kahler
forms) are symplectic deformations of each other, they will (from the \Kahler
point of view) have the same kind of properties in common as if they were
algebraic deformations of each other.

More specifically, this paper represents an extension of the author's
previous results [\WKc ] concerning the deformation properties of the cone of
\Kahler classes $\calK \subset \HZ 2 X$ under deformations of the complex
structure (it may also be viewed as an extension of the work in [\Ru ]). Perhaps
surprisingly, these results generalize in an almost unchanged form to the almost
complex case.  This involves first defining the \Kahler cone for the almost
complex case in terms of pseudo-holomorphic rational curves, and then
generalizing our previous theorem to include almost complex deformations of
$(X,J)$.

For an almost complex
structure $J'$ on $X$, we consider the 
open cone $\HR 2 X$ consisting of classes $D$ with $D^3 >0$ and
such that $D\cdot A >0$ for all homology classes $A$ represented by 
pseudo-holomorphic rational curves   $f : \C \P ^1 \to X $.  We 
\it define \rm the \Kahler cone $\calK \subset \HR 2 X$ of $(X, J')$ to be the
connected component of this cone which contains the classes of symplectic forms 
taming $J'$ (which we note form a convex set).  By results of
[\WPic , \WKc ], this corresponds to the usual \Kahler cone when $(X,J')$ is a
\3f with integrable complex structure.

We shall say that a conic bundle $E$ over a smooth curve $C$ is a \it quasi-ruled
\rm surface over $C$ if the normalization $\tilde E$ is a $\P ^1$-bundle over an
\it unramified \rm cover $\tilde C$ of $C$ -- thus $E$ is either a
ruled surface over $C$ or it is a conic bundle all of whose fibres are line
pairs.  The main result from [\WKc ] said that the \Kahler cone jumped up under
generic complex deformations if the \3f contained quasi-ruled
surfaces over elliptic curves, and otherwise was invariant.  Adopting
terminology from [\WEll ], the \Kahler cone of a generic \it holomorphic \rm
deformation of a \CY threefold $(X,J)$ will be called the \it general
\Kahler cone \rm of $(X,J)$,
and it will equal the \Kahler cone unless $(X,J)$ contains quasi-ruled surfaces
over an elliptic curve.  When $(X,J)$ does contain elliptic quasi-ruled
surfaces, the  \Kahler cone is the subcone of this general
\Kahler cone cut out by the further inequalities $D\cdot l >0$ for $l$ any fibre
of one of the elliptic quasi-ruled surfaces; in [\WEll ], this is proved when
$(X,J)$ contains just elliptic ruled surfaces, but the proof generalizes easily.

\proclaim Theorem 1.  Suppose that $(X,J)$ is a \CY threefold containing 
no quasi-ruled surface over an elliptic curve.  Let
$\omega$ be a \Kahler form on $(X,J)$ and
$\calJ (\omega )$ denote the space of almost complex structures tamed by
$\omega $.
The \Kahler cone of $(X,J')$ at any point $J' \in \calJ (\omega )$ is
a subcone of the \Kahler cone of $(X,J)$.\par\smallskip

\noindent \it Remarks. (i)\quad \rm Since any given \Kahler form $\omega
'$ on $(X,J)$ will tame almost complex structures $J' \in \calJ (\omega )$
sufficiently close to $J$, it follows that the \Kahler cone $\calK$ at $J$ is a
limit of the \Kahler cones of $(X,J')$ as $J'$ tends to $J$.  The \Kahler cone
therefore depends continuously on the almost complex structure at any
integrable point $J$ with $(X,J)$ containing no elliptic quasi-ruled surfaces.
\hfill\par\noindent\it (ii) \quad  \rm The reason why we cannot make the stronger
claim that  $\calK$ equals the \Kahler cone at the generic point $J' \in \calJ
(\omega )$ is that it is at least theoretically possible that there exists a
homology class $A$ (with zero Gromov-Witten invariants) which is $J'$-effective
for $J'$ in the closure of some non-empty open subset of $\calJ (\omega )$, but 
in some open neighbourhood of $J$ 
is not $J'$-effective.  Here we are using the terminolgy from
[\McDS ], where $A$ is called $J'$-effective if it can be written as a positive
integral combination of classes which can be represented by 
pseudo-holomorphic \it rational \rm curves, and we are
implicitly using Gromov
compactness, as stated for instance in [\Sik ], Theorem 5.2.3 or [\McDS ],
Theorem 4.4.3.\hfill\medskip

We can however generalize Theorem 1 to allow deformations of the symplectic
form.  Let $\omega _t$ be any family of symplectic forms on $X$ with $\omega _0$
a \Kahler form on $(X,J)$.  We can then find a family of almost complex
structures $J_t $ on $X$ with $J_0 =J$ such that $J_t$ is always $\omega _t$
tamed.  If $(X,J)$ contains no elliptic quasi-ruled surfaces,
the Main Claim from \S 1 implies that the \Kahler cone of
$(X,J_1 )$ is a subcone of the \Kahler cone $\calK$ of $(X,J)$.  

\proclaim Theorem 2.  Two \CY threefolds which
are deformations of each other as symplectic manifolds will have the same
general \Kahler cone.\par

This follows immediately, since either general \Kahler cone is a
subcone of the other.  There is however an example due to Mark Gross of two
diffeomorphic \CY threefolds (both generic in complex moduli) for which one
\Kahler cone is a proper subset of the other, and hence are not symplectic
deformations of each other [\GrCan , \Ruan ].  Thus each cohomology class which
lies in both \Kahler cones must intersect at least two connected components of
symplectic forms, and the \Kahler forms of the two threefolds must lie in
different components.\smallskip

We can furthermore identify the classes of symplectic forms from the same
connected component as the \Kahler forms on $(X,J)$.

\proclaim Theorem 3.  Assuming that the \3f $(X,J)$ contains no elliptic
quasi-ruled surfaces, the \Kahler cone consists precisely of the classes of
symplectic forms in the same connected component as the \Kahler forms.\par

If now we consider the almost complex structures on $X$ which are tameable
(by some non-specified symplectic form) and which give rise to trivial first
chern class, we obtain a space $\calJ _0$ containing the complex moduli spaces
of all integrable \CY structures on $X$.  Our results show that two such
integrable \CY structures have surprising similarities if they lie in the same
connected component of $\calJ _0$.  Moreover, since such structures in the same
connected component of $\calJ _0$ will have the same \Kahler cone, we can
choose a polarization common to all of them, and then from standard Hilbert
scheme theory there can be at most finitely many complex families contained in
the given component.  I conjecture that only finitely many connected
components of $\calJ _0$ will contain \Kahler structures, and hence
(equivalently) that there are only finitely many families of \CY threefolds for
each diffeomorphism class. I do not in fact know of any examples where more than
one family is contained in a given connected component of $\calJ _0$, and so the
possibility cannot yet be ruled out that two \CY threefolds which are
symplectic deformations of each other might always lie in the same algebraic
family.

The structure of this paper is that our results are reduced to a single Main
Claim in the first Section, which in turn reduces to showing that certain
Gromov-Witten invariants are non-zero.  It might be added here that the
specific results from Symplectic Geometry that we need appear in the
preprints [\Ruan , \RT ], but that we shall also quote the more expository text
[\McDS ].  The Main Claim is checked immediately in the case of
faces of the \Kahler cone corresponding to small contractions.  The remaining
Sections are devoted to proving the Main Claim in the case of non-small
contractions, for which we shall also need a more precise classification of
birational contractions on \CY threefolds.
\bigskip\bigskip

\centerline{\bf 1. \quad  Reduction to Main Claim}\bigskip

\rm  For $(X,J)$ a complex \CY threefold, we denote by $\Kbar$ the closure of
the \Kahler cone $\calK \subset \HR 2 X$, which can be interpreted as
consisisting of the real divisor classes $D$ which are \it nef , \rm that is
$D\cdot C \ge 0$ for all curves $C$ on $X$.  We let $W^* \subset \HR 2 X$
denote the cubic cone of real divisor classes $D$ with $D^3 =0$.  We recall
from the results of [\WKc ] that $\Kbar$ is locally finite rational polyhedral
away from $W^*$, and that the codimension one faces (not contained in $W^*$) of
$\Kbar$ correspond to primitive birational contractions $\phimapX$ of one of
three types.  Type I is when $\phi$ contracts down a finite number of
(numerically equivalent) curves, all smooth and rational.  Type II is when
$\phi$ contracts down an irreducible surface $E$ to a singular point - here,
$E$ will be a generalized del Pezzo surface (see \S 2 for more details).  
Type III is when
$\phi$ contracts down an irreducible surface $E$ to a curve C of Du Val
singularities, generically either an $c{{\rm A}}_1$ or $c{{\rm A}}_2$
singularity, and $E$ is a conic bundle over $C$ (see \S 3 for more details).
In each case, there are rational curves $l$ contracted by $\phi$,
and for any one of these the face is determined by $\{ D\in \Kbar  : 
D\cdot l=0 \}$.  Thus the \Kahler cone as \it defined \rm in the Introduction
in the almost complex case is a natural generalization from the integrable case.

\proclaim Main Claim.  Suppose $(X,J)$ is a \CY threefold.  For every
codimension one face of the nef cone $\Kbar$ not corresponding to the
contraction of an elliptic quasi-ruled surface to a curve of singularities and
not contained in the cubic cone, we can find a homology class $A \in \HZhom 2 X$
which defines the face and is $J_1$-effective for any $J_1$ tamed by a symplectic
form $\omega _1$ in the same connected component as the \Kahler forms on
$(X,J)$.\par

In the case
when the face does correspond  to the contraction of an elliptic quasi-ruled
surface, it follows from [\WKc ] that there exists a class $A$ represented by a
smooth rational curve on $(X,J)$, but not effective for a general holomorphic
deformation.\par

The Main Claim tells us that if $(X,J)$ is general in complex moduli, the \Kahler
cone cannot get larger under a tamed almost complex deformation, and hence
Theorems 1 and 2 from the Introduction follow.  Theorem 3 also follows
easily; the non-trivial claim here is that for any symplectic
form $\omega _1$ in the same connected component as the \Kahler forms on
$(X,J)$, the class of $\omega _1$ lies in the the \Kahler cone $\calK$
of $(X,J)$.
Let $\omega _t$ be a family of symplectic forms on $X$ with $\omega _0$
a \Kahler form on $(X,J)$ and $\omega _1$ the given symplectic form.  As noted
before, we find a family of almost complex structures $J_t $ on $X$ with $J_0 =J$
such that $J_t$ is always $\omega _t$ tamed. If
$(X,J)$ contains no elliptic quasi-ruled surfaces, the Main Claim 
implies that the \Kahler cone of $(X,J_1 )$ is a subcone of the \Kahler cone
$\calK$ of $(X,J)$; however the class of $\omega _1$ clearly lies in the \Kahler
cone of  $(X,J_1 )$, and so the result follows.  The results stated in the
Introduction all reduce therefore to the Main Claim as stated above.
\smallskip

Let us prove the Main Claim in the case of Type I contractions.  This is the
simplest of the three cases, but will illustrate the general principle that in
most cases we can reduce down to considering local \it holomorphic
\rm deformations.  We shall explain the technical results needed from Symplectic
Geometry reasonably fully in this case, and will refer back to the proof below
when considering the other two cases.

\proclaim Proposition 1.1.  If $\phimapX$ is a primitive Type I contraction,
then for each component $Z$ of the exceptional locus, there exists a
neighbourhood $U$ of $Z$ and a holomorphic deformation of the complex structure
on $U$ under which $Z$ splits up into disjoint $(-1,-1)$-curves.\par

\noindent \it Proof.  \rm The morphism $\phi$ contracts $Z$ to a single cDV
(compound Du Val) singularity. 
Taking a Stein neighbourhood $\bar U_0$ of this singularity, there exists by
the argument on page 679 of [\Fr ] a flat 1-parameter family of isolated
threefold singularities $\bar {\cal U} \to \Delta$, where $\bar U_t$ (for $t\ne
0$) contains precisely $\delta >0$ nodes as its singular locus.  Moreover, we
may assume also that there is a simultaneous resolution $\calU \to \Delta$ with
$U_0$ being a neighbourhood of $Z$ in $X=X_0$, and $U_t$ (for $t\ne 0$)
containing precisely $\delta$ disjoint $(-1,-1)$-curves.  We can however take a
\it good representative \rm of our germ $\bar\calU \to \Delta$ such that on the
boundary the family is differentiably trivial (see Theorem 2.8 of [\Loo ]).  If
we take $\calU ' \to \Delta$ to be the corresponding family of resolutions (a
family of manifolds with boundary, where the boundary is differentiably trivial
over $\Delta$), the Ehresmann Fibration Theorem (for manifolds with boundary)
then applies to show that $\calU ' \to \Delta$ is itself differentiably trivial.
Therefore, the family $\calU \to \Delta$ of resolutions may be taken to be 
differentiably trivial, and may as a result be regarded as a family of
holomorphic deformations of the complex structure on the (fixed) neighbourhood
$U$ of $Z$ in $X$.

\medskip  We now use the first feature of almost complex structures which is
not true in the holomorphic case - that they may be patched together in a 
$C ^{\infty }$ way by a partition of unity argument.  Thus, if we take
neighbourhoods of each connected component of the exceptional locus of a
primitive Type I contraction and deform the complex structure locally so that
each component splits up into disjoint $(-1,-1)$-curves, we can patch 
these local (integrable) structures together with the original complex structure
to yield an almost complex structure $J'$ which is integrable in a neighbourhood
of each $(-1,-1)$-curve appearing.  Each such curve will have homology class
lying in the (extremal) 1-dimensional ray determined by the codimension one face.

Let us fix a hyperplane class $H$ of $(X,J)$ and let $\omega$ denote a
corresponding \Kahler form.
If we let $A$ be the homology class represented by one of these 
$(-1,-1)$-curves of minimal degree with respect to the given \Kahler form
$\omega$, then $A$ is $J$-effective, but irreducibly so, since it is only
represented by irreducible embedded holomorphic rational curves.  From Gromov
compactness, we deduce that $J'$ may be chosen above so that $A$ is not
only $J'$-effective but also $J'$-indecomposable, i.e. $A$ is only represented
by simple \curves on $(X,J')$.  Moreover, we may assume that the \curves found
above are the only possible ones
representing $A$.  This latter statement also follows from Gromov compactness,
since by considering the compact subset $F$ of $X$ given by the complement of
suitably small open neighbourhoods of the connected components of the
exceptional locus, we may assume that any other \phrc on $(X,J')$ representing
$A$ intersects $F$.  If this were true for all the almost complex structures
$J'$ constructed in the above way, we would deduce from Gromov compactness that
there was a holomorphic curve on $(X,J)$ representing $A$ and intersecting $F$,
which is nonsense.

\proclaim Proposition 1.2.  The Main Claim is true for primitive Type I
contractions.\par

\noindent \it Proof.\quad  \rm Choose $A$ and $J'$ as above, with $J'$
sufficiently near $J$ that it too is $\omega$-tamed (this being an open
condition).  A standard regularity criterion (see [\McD  ] (4.2) or [\McDS ]
(3.5.1)) then implies that the almost complex structure $J'$ is regular (for the
class $A$) -- note, that the proof given only requires $J'$ to be integrable in a
neighbourhood of each \phrc representing $A$.

Since $(X,J)$ is Calabi--Yau, 
the symplectic manifold $(X, \omega )$ is weakly monotone, and any almost
complex structure in $\calJ (X, \omega )$ is semi-positive.  Thus if $A$ is a
non-multiple class, the theory of [\McDS ], \S 7.2, gives a well-defined 
3-point Gromov-Witten invariant $\Phi _A (H,H,H)$, invariant also under
symplectic deformations of $\omega $.
Moreover, the property that $J'$ is integrable in a
neighbourhood of each \phrc representing $A$
implies that $J'$ induces an (integrable) complex structure on the moduli space
$\calM (A, J')$ (see [\McDS ], Remark 3.3.6), and hence the natural orientation
at the points in the moduli space of unparametrized curves representing $A$ is
always positive (cf [\McDS ], Remark 7.3.6).  
Thus $ \Phi _A (H,H,H) = \Phi _{A,J'} (H,H,H)$ is equal to a positive
multiple of $(A\cdot H)^3$. 

When $A$ is a multiple class, i.e. $A=mB$ for some $m>1$, the Gromov
-Witten invariant  $ \Phi _A  $ is not well-defined for the reasons explained
in [\McDS ], \S 9.1.  To circumvent this difficulty, we follow the
original idea of Gromov and work with the graphs of \curves $\C\P ^1 \to X$.
We let $\hat X$ denote the product manifold $\C \P ^1 \times X$, with $\hat
\omega$ the product symplectic form $\tau _0 \times \omega$ where $\tau _0$ is
the standard symplectic form on $\C \P ^1$ corresponding to the Fubini-Study
metric, and $\hat J'$ the product almost complex structure $i \times J'$ with
$i$ the complex structure on $\C \P ^1$.  Let $A_0$ denote the homology class 
$[\C\P^1 \times \{ pt \} ] \in \HZhom 2 {\hat X}$, and let $\hat A$ denote the
class $ A_0 + A \in \HZhom 2 {\hat X}$.  As explained in [\McDS ], \S 9.1, we now
work on $\hat X$ rather than $X$, counting $\hat A$-curves on $\hat X$ rather
than $A$-curves on $X$; the advantage of this is that the class $\hat A$ is
necessarily primitive (i.e. non-multiple), whilst the same was not true for $A$. 
For a general $J'' \in \calJ (\omega )$, we must then
work, not with the product almost complex structure, but with a small  
perturbation.

The theory described in [\McDS ] does not quite cover the \CY
case, and so instead we refer to [\Ruan ].  The invariant we need is called
$\tilde \Phi$ there, and its construction and properties are described full in
\S 3 of the paper -- it is in fact also a special case of the mixed invariant
$\Phi $ constructed in [\RT ], which could therefore serve as an alternative
reference.  To define $\tilde \Phi$, Ruan follows the idea of Gromov and
considers a special type of small perturbation $J_g$ from the product almost
complex stucture on $\C \P ^1 \times X$, where $g$ is an anti-complex linear 
bundle map from $T\C\P ^1$ to $TX$.   
With this structure, a map $f : \C\P ^1 \to X$ satisfies the equation 
$\bar \partial _{J''} f = g$ (which should be regarded as a small
perturbation from the condition $\bar \partial _{J''} f = 0$ needed for $f$ to be
holomorphic) if and only if the corresponding
section $\bar f : \C \P ^1 \to \C \P ^1 \times X$ is  $J_g$-holomorphic.  

In this way,
he is able to construct a well-defined invariant $\tilde \Phi _{(A, \omega )}$
by ``counting the number of perturbed holomorphic maps with marked points" (see
[\Ruan ], (3.3.6)), which in our case will also be invariant under symplectic
deformations of $\omega $ (see Lemma 3.3.4 and note that semi-positivity is
satisfied in our case).  The fact that the invariant is well-defined essentially
follows from the fact that $\hat A$ is a primitive (i.e. non-multiple) class and
a standard cobordism argument. We shall in this paper therefore denote this
invariant by $\tilde \Phi _A$, and refer to $\tilde \Phi _A (H,H,H)$ as the
three point Gromov-Witten invariant -- it might be added 
however that the literature in
this field employs many different terminologies and notations.  In
the case when $A$ is primitive, it follows from Proposition 3.3.8 of [\Ruan ]
that two definitions yield the same three point invariant.

In order to calculate this invariant in our particular case, we observe that the
pair $(J' ,0)$ is $A$-good, in the terminology of [\Ruan ], (3.3.7) -- this
follows from the $A$-regularity of $J'$ observed before.  Thus we are able to
calculate $\tilde \Phi _A (H,H,H)$ from the almost complex structure $J'$.  In
fact, we have a precise knowledge of the moduli space $\calM _{(A,J',0)}$ of
parametrized \curves representing $A$, and if we consider as in [\Ruan ] marked
\curves  (with three marked points on $\C\P^1$) then we obtain a finite set. 
This enables us to calculate that $\tilde \Phi _A (H,H,H) = \tilde \Phi
_{(A,J',0)}  (H,H,H)$ is equal to some positive multiple of $(A\cdot H)^3$
as before.

Thus, for small perturbations $J _{1,g}$ of the product almost
complex structure $i \times J_1 $, there are \curves representing the class
$\hat A = A_0 + A$.  By Gromov compactness, we deduce for the product
structure that the the class $\hat A$ is represented by a \phrc or cusp curve,
and hence projecting onto the first factor that the class A is $J_1$-effective. 
The proof of the Main Claim is now complete for Type I contractions.

The above argument provides the model also for the proof of the Main
Claim in the case of contractions of Types II and III; we find a class $A$
defining the relevant face of $\Kbar$ for which $\tilde \Phi _A (H,H,H)$ can be
shown non-zero, and hence the Main Claim will follow.\bigskip\bigskip

\centerline{\bf 2. \quad  Type II contractions}\bigskip

\rm  Suppose now that $\phimapX$ is a primitive Type II contraction, therefore
contracting a generalized del Pezzo surface $E$ down to a singular point.  As
was the case for Type I contractions, the technique we use to prove the Main
Claim is to make a small holomorphic deformation on some open neighbourhood of
the exceptional locus, under which $E$ deforms to a smooth del Pezzo surface --
the fact that this is always possible has been proved by Gross in [\GrDef ].  We
then patch together with the original complex structure, and observe that in
most cases there is an obvious choice for the class $A$.  Where this is not the
case, we have to use a cobordism argument of Ruan.

In passing, I remark that it is an interesting question whether there is a
global holomorphic deformation of $X$ under which $E$ deforms to a smooth del
Pezzo surface.  This is unknown and is the analogous problem to a conjecture of
Clemens for Type I contractions that there is a global holomorphic deformation
under which the exceptional locus splits up into 
distinct $(-1,-1)$-curves -- cf. Proposition 1.1
where we prove that this does happen for deformations on some open neighbourhood
of the exceptional locus.

Generalized del Pezzo surfaces may be classified, but as we see below, not all
of them occur as the exceptional surface of a Type II contraction.  The normal
del Pezzo surfaces are either elliptic cones or del Pezzo surfaces with \rdp
singularities as classified in [\Dem  ] -- for a proof of this statement, see for
instance [\HW ].  The non-normal del Pezzo surfaces have been classified in
[\ReDP ]. Let us first consider the possibilities for small values of $E^3$.

If $E^3 =1$, the singularity of $\Xbar$ is a hypersurface singularity with
equation of the form $x^2 +y^3 + f(y,z,t) =0$ where $f= yf_1 (z,t) +f_2 (z,t)$
and $f_1$ (respectively $f_2$) is a sum of monomials of degree at least 4 
(respectively 6).  The morphism $\phimapX$ is then the weighted blow up of
the singularity with weighting $\alpha$ given by $\alpha (x) =3$, 
$\alpha (y) =2$ and $\alpha (z) = \alpha (t) =1$ (for a proof of these
statements, see [\RCan ] (2.10) and (2.11)).  Thus $E$ is naturally embedded in 
$\P (3,2,1,1)$ with equation $x^2 +y^3 + yg_1 (z,t) +g_2 (z,t)=0$, where $g_1$ 
(resp. $g_2$) is the degree 4 (resp. 6) homogeneous part of $f_1$ (resp.
$f_2$).  It is clear however that the hypersurface singularity may be deformed
locally so that the exceptional locus $E_1$ of the deformed singularity is a
smooth del Pezzo surface of degree one in $\P (3,2,1,1)$.  Note that the
singularities can be resolved in the family by means of the weighted blow up
described above.

By taking a good representative of this deformation and applying Ehresmann's
fibration theorem (with boundary), as was done in the proof of (1.1), we can
find a neighbourhood $U$ of $E$ in $X$ and a holomorphic deformation of the
complex structure on $U$ such that $E$ deforms to a smooth del Pezzo surface of
degree one.

The case $E^3 =2$ is precisely analogous, with $\Xbar$ having a hypersurface
singularity with equation $x^2 + f(y,z,t) =0$ with $f$ a sum of monomials of
degree $\ge 4$, and the desingularization given by the weighted blow up with
weighting $\alpha $ given by $\alpha (x) =2$, $\alpha (y) =\alpha (z) =
\alpha (t) =1$ (again see [\RCan ], (2.10) and (2.11)).  Thus $E$ is naturally
embedded in $\P (2,1,1,1)$ with equation $ x^2 + f_4 (y,z,t) =0$.  The previous
argument applies again to show that by deforming the complex structure on a
neighbourhood $U$ of $E$, we can deform $E$ to a smooth del Pezzo surface of
degree two.

For $k =E^3 \ge 3$, we know that the morphism $\phimapX$ is an ordinary blow
up of the singular point, 
and $\Xbar$ has local embedding dimension $k+1$ [\RCan ].  In the case of $k=3$ for
instance, $\Xbar$ still has a hypersurface singularity and $E$ is a cubic
hypersurface in $\P ^3$.  Again therefore, we can deform the complex structure
on a neighbourhood $U$ of $E$ so that $E$ deforms to a smooth cubic surface.

\proclaim Lemma 2.1.  Suppose $\phimapX$ is a primitive Type II
contraction of a surface $E$ with $k=E^3 \le 3$; then the Main Claim from \S 1
is true for the corresponding face of the nef cone.\par

\noindent \it Proof. \rm  We can find a neighbourhood $U$ of $E$ and a
holomorphic deformation of the complex structure on $U$ such that $E$ deforms
to a smooth del Pezzo surface.  Such a surface however has a finite number of
$(-1)$-curves, and these are $(-1,-1)$-curves in the deformed complex structure
on $U$.  Let $A$ be the homology class of one of these curves, a primitive
class since $E\cdot A = -1$.

By patching this deformed complex structure on $U$ together with the original
one on $X$ (as was done in (1.2)), we obtain an almost complex structure $J'$
on $X$.  Furthermore, we may assume that the above $(-1,-1)$-curves are the
only \curves on $(X,J')$ representing $A$ (the same argument via Gromov
compactness as in the preamble to (1.2)), and that $J'$ is regular for the
class $A$.  The proof is now similar to (1.2), showing that the 3-point
Gromov-Witten invariant $\Phi _A (H,H,H)$ is positive and that $A$ is
$J_1$-effective for any $J_1$ as in the statement of the Main Claim.  More
precisely, the proof is a simpler version of (1.2) by reason of the
fact that the homology class $A$ is primitive, and so the
Gromov-Witten invariant $\Phi _A (H,H,H)$ may be used.\medskip

\noindent \it Example 2.2. \rm   Consider the case when the exceptional
surface $E$ of the contraction $\phimapX$ is a projective cone on a smooth plane
cubic.  We can then deform the complex structure locally in a neighbourhood of
$E$ so that $E$ deforms to a smooth cubic surface, containing 27 lines.  The
above argument shows that these 27 \curves persist under almost complex
deformations of the complex structure on $X$.  A classical analysis of the
degeneration of the 27 lines from the smooth cubic (see for instance [\BW ])
reveals that they degenerate three at a time to the 9 generators of the cone $E$
corresponding to the inflexion points of the plane cubic.  Thus it is these
generators which deform under an almost complex deformation of the structure on
$X$, generically each deforming to three pseudo-holomorphic rational
curves.\bigskip

\proclaim Lemma 2.3.  Suppose $(V,P)$ is an isolated rational Gorenstein 3-fold
singularity for which blowing up the point $P$ yields a crepant resolution 
$(U,E)$, with the exceptional divisor $E$ irreducible and $E^3 =k >3$.  Then $E$
is either a normal del Pezzo surface of degree $k\le 9$ with at worst rational
double point singularities, or a  non-normal del Pezzo surface of degree $k=7$
whose normalization is a non-singular rational ruled surface.\par

\noindent \it Proof.\rm \quad  We can reduce immediately to the case when $E$ is 
non-normal; if $E$ is normal, it cannot be an elliptic cone, since then $E$ 
would have embedding dimension $>3$ at its vertex, contradicting the assumption
that $U$ is smooth.   The only remaining possibility in the normal case is that
$E$ is a del Pezzo surface of  degree $\le 9$ with only rational double point
singularities.

Suppose then that $E$ is non-normal with $E^3 =k>3$.  Turning to the 
classification of non-normal del Pezzo surfaces
in Theorem 1.1 of [\ReDP ], we observe that case (c) (projection of a cone $\F
_{k;0}$) does not occur, since the point of $E$ corresponding to the vertex of 
$\F _{k;0}$ would have too high an embedding dimension.  We are left therefore
with the cases when $E$ is either a projection $\bar \F _{a;1}$ of a scroll 
$\F _{a;1}$ ($k=a+2$ where we assume $a\ge 2$), or the 
projection $\bar \F _{a;2}$ of a scroll 
$\F _{a;2}$ ($k=a+4$ where we assume $a\ge 0$).

In the first case the projection of  $\F _{a;1}$ is taken with centre a general
point in the plane spanned by the minimal section and a fibre (both being
lines), and the resulting surface $\bar \F _{a;1} \subset \P ^{a+2}$ has a
double line.  In the second case, 
the projection of  $\F _{a;2}$ is taken with centre a general
point in the plane spanned by the minimal section (a conic).  The resulting 
surface $\bar \F _{a;2} \subset \P ^{a+4}$ again has a double line.

Let us consider the case first when $E \cong \bar \F _{a;2}$.  It was observed
first by Gross in [\GrDef ] that this forces $a=3$ and hence $k=7$; we give
here a different proof (part of which will be needed in \S 4) of this fact.

Let $C$ denote the line of singularities of $E$ and suppose its normal bundle
is $N_{C/U} = (c,-2-c)$, with the usual notation that $(r,s) = \calO _{\P^1}
(r) \oplus \calO _{\P^1} (s)$.  If $g : \tilde U \to U$ denotes the blow up of
$U$ in $C$, we obtain on $\tilde U$ a scroll $E' \cong \F _e $ where $e=2c+2$. 
If $A$ denotes the minimal section of $E'$, we have $K_{E'} \sim -2A - (e+2)f$
(where $f$ denotes a fibre), and hence  $-E' |_{E'} \sim A + (c+2)f$.  Moreover,
if $\tE$ denotes the strict transform of $E$, it is easily checked that the
morphism $\tE \to E$ coincides with the projection $ \F _{a;2} \to \bar \F
_{a;2}$, and hence that $\tC = \tE \cap E'$ is both isomorphic to $\P ^1$ and is
a double section of $E' \cong \F _e$. From this, an easy calculation shows that
$\tC \sim 2A + (e+1)f$.  But $\tC$ is irreducible and so $\tC \cdot A \ge 0$,
from which we deduce $e\le 1$, and hence (since $e$ is even) that $e=0$.  Thus
$c=-1$ and the normal bundle $N_{C/U} = (-1,-1)$.

To show that $a=3$ and hence $k=7$, we now calculate on $\tE \cong \F _{a;2}$. 
Note first that $(\tC ^2 )_{E'} = 4A^2 + 4(e+1) = 4$.  Now 
$$ -2 = \deg K_{\tC } = (K_{\tE } + \tC )\cdot \tC = (\tE +2E')\cdot \tC
= 4 + 2 (\tC ^2 )_{\tE }.$$
Thus $ 2 (\tC ^2 )_{\tE } = -6$, i.e.~$(\tC ^2 )_{\tE } = -3$.  
Hence $a=3$ as claimed.

A similar elementary calculation shows that only $E^3 =7$ could occur for the
case $E \cong \bar \F _{a;1}$ with $a\ge 2$, but this in fact also follows from
more general considerations.  If the blow-up $(U,E) \to (V,P)$ is as in the
statement of the Proposition (without any restriction on $k$), it is shown in
the proof of Theorem 5.8 of [\GrDef ] that any deformation $E_1$ of $E$ may be
achieved by deforming the complex structure on $U$, which may in turn be
achieved by deforming the singularity $(V,P)$.  But for $a\ge 0$ the surface
$\bar \F _{a+2;1} $ in $\P ^{a+4}$ is a degeneration of surfaces isomorphic to 
$\bar \F _{a;2} $.  To see this, we consider the scroll $\F _{a+1;2} 
\subset \P ^{a+6}$ given parametrically by $(\lambda s^{a+3} : \lambda s^{a+2} t
: \ldots : \lambda t^{a+3} : \mu s^2 : \mu st : \mu t^2 )$, for $(\lambda :\mu
) \in \P^1$ and $(s:t) \in \P^1$, and project from points $(0:\ldots :0:
u:0:0:v )$ on the ruling given by $s=0$
to get a family of surfaces in $\P ^{a+5}$.  If we now
project this family from the point $(0:\ldots :0;-1;1;1) \in \P ^{a+5}$, we get
a family in $\P ^{a+4}$ whose general member is a $\bar \F _{a;2} $, but which
degenerates to a $\bar \F _{a+2;1} $ (corresponding to the first projection
having centre the point $(0:\ldots :0:1)$ on $\F _{a+1;2} $). As remarked
in [\GrDef ], this fact together with the previous result shows that only the
case $E \cong \bar \F _{5;1}$ can occur.\bigskip

\proclaim Proposition 2.4.  If $\phimapX$ is a primitive Type II
contraction on a \CY threefold $X$, then the Main Claim from \S 1
is true for the corresponding face of the nef cone.\par

\noindent \it Proof.\rm\quad  By (2.1), we may assume that $k= E^3 >3$.  Suppose first 
that $E$ is non-normal, and so $E \cong \bar \F _{3,2}$ or $\bar \F _{5,1}$ in the
notation of (2.3).  It is easy to see ([\GrDef ], (5.6)(iii)) that $\bar \F
_{3,2}$ can be smoothed in $\P ^7$ to a smooth del Pezzo surface $E_1$ of degree
7, and hence by the argument in the proof of (2.3), the same is true for $\bar
\F _{5,1}$.  Moreover, we may  apply again the argument from Theorem 5.8 of
[\GrDef ] to deduce that such a smoothing  $E_1$ of $E$ may be achieved by
holomorphically deforming the complex structure  on some neighbourhood $U$ of
$E$.  The surface $E_1$ is a smooth del Pezzo  surface of degree 7 containing
precisely three $(-1,-1)$-curves.   (We remark that in the case when  $E \cong
\bar \F _{3,2}$, one of these $(-1,-1)$-curves corresponds to the singular 
locus of $E$, which we saw in (2.3) was already a $(-1,-1)$-curve, and the
other  two correspond to the double fibres of $\bar \F _{3,2}$.)    We let $A$
be the homology class of such a $(-1,-1)$-curve.   We may patch together the
deformed holomorphic complex  structure on $U$ with the original complex
structure $J$ to get an almost  complex structure $J'$ on $X$, exactly as was
done in (2.1).   The proof of (2.1) now goes over unchanged to
prove the Main Claim in this case.

By (2.3), we have reduced to the case when $E$ is a del Pezzo surface with at 
worst rational double point singularities.  Let us
dispose first of the cases when $E$ is neither $\P^2$, nor $\P^1 \times \P^1$,
nor the quadric cone in $\P ^3$. Here we can proceed in one of two ways.  On the
one hand, we can observe that $E$ will contain lines with normal bundle
$(-1,-1)$ or $(-2,0)$ in the threefold.  In the latter case, the curves do not
move in the threefold and so may be locally analytically contracted [\RMin ]. 
The argument from \S 1 then shows that such curves split up into disjoint
$(-1,-1)$-curves when we make a local deformation of the complex structure, and
the Main Claim follows as in (1.2).  The alternative argument is the one we have
just used above, where we observe  that $E$ may be smoothed, and that such a
smoothing  may be achieved by holomorphically deforming the complex structure on
some  neighbourhood $U$ of $E$.  Under such a deformation, $E$ deforms to a
smooth  del Pezzo
surface $E_1$, not isomorphic to either $\P^2$ or $\P^1 \times \P^1$, and 
hence containing $(-1,-1)$-curves.  The proof from (2.1) now goes over
unchanged to prove the Main Claim in these cases.

We consider now the cases when $E$ is 
isomorphic to either $\P^2$ or $\P^1 \times \P^1$.  In both cases, we have
smooth families of $\P ^1$s with no distinguished members, and we need to know
how many of these curves will persist is we make a generic deformation of the
almost complex structure on $X$.  There is an intuitively attractive argument
in \S 8 of [\CFKM ] that this number is $(-1)^b e(B)$, where $B$ is the
smooth Hilbert scheme parametrizing the family, $e(B)$ the Euler
characteristic, and $b$ the dimension of $B$.  The point is that the
obstruction bundle to first order deformations is, as commented in [\CFKM  ],
just the cotangent bundle to $B$, and a $C^{\infty}$ 1-parameter deformation
gives rise to a $C^{\infty}$ 1-form on $B$; the above number is just the
expected number of zeros of such a 1-form (taking into account orientations). 
The problem with this approach is that, unlike the algebraic case, it will not
be enough to show that a given member of the family deforms to all finite
orders; we will also have convergence to worry about.  The way round this is to
invoke a cobordism argument due to Ruan.

Consider the case when $E$ is isomorphic to $\P ^2$.  We take $A$ to be the
homology class of a line on $E$ -- note that $E \cdot A =-3$.  We now refer the
reader to Proposition 5.7 of [\Ru ].  Note that although $A$ is not a primitive
class in homology, it is both extremal in the cone of 1-cycles and not a
non-trivial multiple of any $J$-effective class.  Arguing from Gromov
compactness, the same will be true for almost complex structures in some
neighbourhood of $J$, and so for almost complex structures in such a
neighbourhood we do not need to worry about multiple 
pseudo-holomorphic rational curves.  The result from [\Ru ] does apply
therefore to our case, even though $A$ is not primitive.  The other conditions
of the Proposition are now easily checked to be satisfied, recalling that 
for $f : \C \P ^1 \to X$ a holomorphic rational curve, 
$\dim\, coker (\bar D ) = h^1 (f^* T_X)$,
and so the required equality holds because the Hilbert scheme $B$ is smooth (in
this case just $\P ^2$).  We deduce that for any generic compatible almost
complex structure $J'$ sufficiently close to $J$, the moduli space ${{\cal M }}
_{A,J'}$ of non-parametrised \curves is oriented cobordant to the zero set of a
tranverse section of the cotangent bundle to $B$, and so consists of $e(B)=3$
points when counted taking into account orientations.  Thus for $H$ a
hyperplane class on $(X,J)$, it follows as in (1.2) that the Gromov-Witten
invariant $\tilde \Phi _A (H,H,H)$ is strictly positive, and so the Main Claim
follows.

The case of $E$ isomorphic to $\P^1 \times \P^1$ is entirely analogous.  If $A$
denotes the class of a line on the quadric $\P^1 \times \P^1$, the
corresponding Hilbert scheme consists of two disjoint lines, and in this case
the Gromov-Witten invariant $\tilde \Phi _A (H,H,H)$ will be strictly negative.  
The remaining case, when $E$ is isomorphic to the quadric cone, may be reduced 
to this case, if we smooth $E$ by means of a local holomorphic deformation on
some neighbourhood of $E$ and patch together with the original complex structure
$J$ on $X$.  The Main Claim is therefore proved for all Type II contractions.
\bigskip\bigskip

\centerline{\bf 3. \quad Classification of Type III contractions}\bigskip

\rm  
\proclaim Proposition 3.1.  If $\phimapX$ is a primitive Type III contraction
on a smooth \CY threefold $X$, then the exceptional locus $E$ is a conic bundle
over a smooth curve $C$, where $\phi$ is just contraction along fibres on
$E$.\par

\noindent \it Proof. \rm  This is a natural generalization of an argument from
pages 568-9 of [\WKc ], and we shall adopt the same notation as used there. 
In particular $E$ will denote the exceptional divisor on $X$ with $C$ its image
on $\Xbar $; set $\pi : E \to C$ to be the induced morphism.  Letting $D$ denote
the  pullback of a hyperplane section of $\Xbar $, we have for $n$ sufficiently
large that $nD-E$ is ample and that the map $H^0 (X, \calO _X (nD)) \to 
H^0 (E, \calO _E (nD))$ is surjective, from which it follows that 
$\pi _* \calO _E = \calO _C$.
Under our assumptions, it
was observed in [\WKc ] that the sheaf
$\calO _X  (-E)$ is relatively generated by its global sections, i.e.~that the
map 
$\phi ^* \phi _* \calO _X (-E) \to \calO _X (-E)$ is surjective.  From this, it
followed that the map $\phi $ has to factor through the blow-up $X_1$ of $X$ in
$C$, and hence (from the primitivity of $\phi $) if $X_1$ is normal, we may
take $\phi $ to be the blow-up morphism.  Furthermore the linear system $|nD-E|$
is without fixed points for $n$ sufficiently large.  

Once we have shown that
$C$ is smooth, an easy argument (cf. calculations below) then shows that $X_1$ is
locally a complete intersection and non-singular in codimension one, and
therefore normal; hence
$X$ is the blow-up of $\Xbar $ in $C$. Since the intersection number of $-E$
with the general fibre of $E$ over $C$ is two, and all curves contained in
fibres are in the same numerical ray, it follows (from the fact that $C$ is a
curve of cDV singularities) that $E$ is a conic bundle over $C$.  From this we
see that the divisor $nD-E$ will in fact be very ample for $n$ sufficiently
large.

To show that $C$ is smooth, we suppose that $P$ is a point of $C$ and $Z$ is the
corresponding reduced fibre over $P$; then $Z$ consists of one or two components
(in all cases isomorphic to $\P ^1$).  Taking $S$ to be a general
element of $|nD -E|$ for $n$ sufficiently large, we observe that $S$ 
can intersect $Z$ in one or two points.  The case when $S$ intersects
$Z$ in two points $Q_1$ and $Q_2$ has been dealt with in [\WKc ].  For this
case, we observe that $S$ is smooth at both points, and  
the projection $\pi : E \to C$ yields a local analytic isomorphism
between the curve $E\cap S$ at $Q_i$ and the curve $C$ at $P$.  Thus we deduce
that if $P$ were a singular point of $C$, it would only have embedding dimension
two.  It is however observed on p.~569 of [\WKc ] that the blow-up $X_1$ of 
$\Xbar$ in $C$ is non-singular in
codimension one and locally a complete intersection (by an argument
given below, it is locally principal on $\tilde \C ^4 $), and hence normal, from
which it follows that the blow-up must in fact be $X$.  The calculation
performed on p.~569 of [\WKc ] shows however (under the assumption that $P$ has
embedding dimension two) that this blow-up fails to be smooth when $P$ is a
singular point of $C$; we deduce therefore that $C$ is smooth.

We are left with the case when the 
general element $S$ of $|nD-E|$ intersects $Z$ in one point, where $Z$ is now
isomorphic to $\P ^1$.  Setting $\tC = S\cap E$, we know that 
$\tC$ is a double cover of $C$ with only a single point $Q$ lying above
$P$.  The embedding dimension at $Q$ of $\tC$ is at most two, and there is an
involution $\iota : \tC \to \tC$ switching the two points of a general fibre.
If moreover $Q$ is a non-singular point of $\tC$, then $P$ is a non-singular
point of $C$;  we shall assume from now on that $Q$ is singular.
Therefore locally analytically we have an induced involution $\iota : \C ^2
\to \C ^2$ with $\tC$ defined by an invariant (or anti-invariant) element $h \in 
\C [[\eta ,\xi ]]$.

The involution $\iota$ on $\C ^2$ may be diagonalized so as to act via the
matrix $\bigl({1 \atop 0}{0 \atop -1}\bigr)$ or $\bigl({-1 \atop 0}{0 \atop
-1}\bigr)$; in the former case, the quotient is still $\C ^2$ and so $C$ has
embedding dimension two again and the previous argument from [\WKc ] holds.
We should therefore concentrate on the latter case.

Let us assume first that $\tC$ is defined by an invariant element $h \in 
\C [[\eta ,\xi ]]$, whose
terms are therefore all of even degree. 
We know that $P\in \Xbar$ has embedding dimension 4; 
choose local coordinates $x,y,z,t$ on $\C ^4$ so that the quotient of $S$ by
the the induced involution $\iota$ has image given by $x=0$, $yz=t^2$
(i.e.~where $x=0$ defines the smooth surface $S$ on $X$ and $t=\eta \xi$,
$y=\eta ^2$, $z=\xi ^2$).  The function $h$ may now be considered as an element
of $\C [[y,z,t]]$, and $C$ is cut out by $h$ on the surface given by $x=0$, $yz
=t^2$ in $\C ^4$.  The
blow-up $\tilde \C ^4$ of $\C ^4$ in $C$ may be considered (analytically) as
the subvariety of $\C ^4 \times \P ^2$ defined by 
$$ \rank \pmatrix{x & yz-t^2 & h \cr
                  u &    v   & w \cr} \le 1,$$
where $u,v,w$ are
homogeneous coordinates on $\P ^2$.  Taking the affine piece of $\tilde \C ^4$
given by $u=1$, we get equations for this affine piece as a subvariety of $\C
^6 = \C ^4 \times \C ^2$ to be $yz-t^2 = xv$, $h=xw$.  Calculating the partial
derivatives of the two equations and evaluating them at any point above $P =
(0,0,0,0) \in \C ^4$, we get row vectors $(v,0,0,0,0,0)$ and $(w,-a,-b,-c,0,0)$
where $h=ay+bz+ct+$~{higher order terms}.  Let $L$ denote the line in the
locus $P\times \P ^2 \subset \tilde\C ^4$ given by $v=0$; then the
above calculation shows that $\tilde \C ^4$ is singular at all points of $L$. 
Observe that $\tilde\C ^4$ is non-singular in codimension one and
locally a complete intersection, hence normal.  Since the exceptional divisor
on $\tilde \C ^4$ is locally principal (from the fundamental property of
blow-ups), it follows that the strict transform $X _1 \subset \tilde \C ^4$ of
$\Xbar \subset \C ^4$ is locally principal, where $X_1$ is of course also the
blow-up of $\Xbar $ in $C$.  Thus $X_1$ is itself a local complete
intersection.  As argued before however, the map
$\phimapX$ factors through this blow-up, and so $X_1$ is non-singular in
codimension one and in particular normal.  From the primitivity of $\phi$, we
deduce that $X$ is isomorphic to $X_1$, i.e.~
$X$ is just the blow-up of $\Xbar $ in the curve $C$. 
Moreover, the exceptional locus of $\phi$ above $P$ is the intersection of $X$
with  $P\times \P ^2 \subset \tilde\C ^4$, and this is just the
fibre of $E$ above $P$, i.e.~in our case the curve $Z$.  Thus $X$ has a
singularity where $L$ meets $Z$ in $\P ^2$, contradicting our assumption that
$X$ was smooth.

The remaining case to consider is when the involution $\iota$ on $\C ^2$
again acts via the matrix $\bigl({-1 \atop 0}{0 \atop -1}\bigr)$ but 
$\tC$ is defined by an anti-invariant element $h \in \C [[\eta ,\xi ]]$, whose
non-zero terms are therefore all of odd degree.  So $C \subset \C ^4$ is
defined by equations $x=0$, $yz=t^2$, $h_1 =0$, $h_2 =0$, where 
$ h_1$, $h_2 \in \C [[y,z,t]]$ and the images of $h_1$, $h_2$ in 
$\C [[\eta ,\xi ]]$ are respectively $\eta h$ and $\xi h$.  In particular, we
have the relations $$ th_1 = yh_2 + (t^2 -yz)g_1 \hbox{\rm \quad and \quad  }
th_2 = zh_1 + (t^2 -yz)g_2$$
for suitable $ g_1$, $g_2 \in \C [[y,z,t]]$.

The blow-up $\tilde \C ^4$ of $\C ^4$ in $C$ may then be
considered as a component of the subvariety of  $\C ^4 \times \P ^3$ defined by 
$$ \rank \pmatrix{x & t^2 -yz & h_1 & h_2 \cr
                  r &    u    & v   & w \cr} \le 1,$$
where $r,u,v,w$ are
homogeneous coordinates on $\P ^3$.  Note that these equations define a
reducible variety, with one component (which we do not want) having fibre $\P ^3$
over all points of $C$.  Observe now that the equations imply that 
$xtv = rth_1 = ryh_2 + r(t^2 -yz)g_1$, where $ryh_2 =xyw$.  Thus $x(tv-yw) = 
r(t^2 -yz)g_1$.  Multiplying by $u$, we deduce that 
$$r(t^2 -yz)(tv - yw -ug_1 ) =0,$$ and hence that $tv-yw -ug_1$ is in the ideal
of definition for $\tilde \C ^4 \subset \C ^4 \times \P ^3$.  Similarly, we
show that $tw - zv - ug_2$ is in the ideal of definition.  Let 
$W \subset \C ^4 \times \P ^3$ be the subvariety defined by  equations 
$xu = r(t^2 -yz)$, $tv - yw - ug_1 =0$, $tw-zv-ug_2 =0$, $xv =h_1$, $xw=rh_2$. 
Note that the fibre of $W$ above a general point of $C$ is $\P ^2$,
whilst above the point $P$ it will be $\P ^3$ when $g_1 (P) = 0 = g_2 (P)$, 
and $\P ^2$ otherwise. 
Observe that if $h = a\eta + b\xi + \hbox{\rm higher order terms}$, then $g_1
(P) = b$ and $g_2 (P) =a$; thus our assumption that $Q$ is singular ensures
that $g_1 (P) = 0 = g_2 (P)$.

Taking partial derivatives of the above set of defining equations, and
evaluating the corresponding Jacobian matrix at any point of $W$
above the origin $P\in \C ^4$, it is easily checked that the corresponding
matrix has rank 3 unless $v=w=0$.  Since the blow-up $\tilde \C ^4$ coincides
with $W$ outside $P \times \P ^3$, it follows that $W = \tilde \C ^4$.
The above calculation has shown that $\tilde \C ^4$ has fibre $\P ^3$ above $P$
and is smooth except possibly along the line $v=w=0$ in $P\times \P ^3$.  We
consider the contraction $\Xbar \subset \C ^4$, a hypersurface containing $C$
and having only cDV singularities along $C$ -- in fact generically along $C$
either 
${{\rm cA}}_1$ or ${{\rm cA}}_2$ singularities.  Such a hypersurface will have
a local analytic equation of the form 
$$ 0 = F(x,y,z,t) = \alpha x^2 + \beta x(t^2-yz) +\gamma (t^2-yz)^2 + \delta _1
xh_1 + \delta _2 xh_2 $$ 
\rightline {$\quad + \sigma _1 (t^2-yz)h_1 + \sigma _2 (t^2-yz)h_2
+ \hbox{\rm terms involving } h_1^2, h_1h_2, h_2^2,$}\medskip

\noindent where $\alpha (P) \ne 0$ since $F$ is not an element of $m_P ^3$ 
(recalling that $P$ is a cDV singularity).

If we look at the affine piece of the blow-up given by $r=1$, the strict
transform of $\Xbar$ (i.e.~the blow-up $X_1$) is given by an equation of the form
$$ 0 = \alpha  + \beta u +\gamma u^2 + \delta _1 v + \delta _2 w  
+ \sigma _1 uv + \sigma _2 uw
+ \hbox{\rm terms involving } v^2, vw, w^2.$$
So $X_1$ is cut out on $\tilde \C ^4$ by a single extra equation, which
therefore cuts out a two dimensional subvariety on $P\times \P ^3$, not
containing the line $v=w=0$ in $P\times \P ^3$.

Under the assumption that $P$ is a singular point of $C$, the above analysis
has shown that the fibre of $X_1$ above $P$ is two-dimensional.  As observed
before, $\phimapX $ factors through
the blow-up $X_1$ of $\Xbar$ in $C$.  Since the fibre of $X$ above $P$ is
1-dimensional, we arrive at the required contradiction.  The proof of (3.1) is
now complete. \bigskip

Having seen that $C$ is smooth, we can make a local calculation to see what are
the possible dissident singularities.  In the case when the generic singularity
is an ${{\rm cA}}_1$ singularity, the restrictions obtained are rather weak (see
[\WKc ], Section 2, for a discussion of this case).  In the case when the 
generic singularity
is a ${{\rm cA}}_2$ singularity, the restrictions on a dissident singularity
$P$ (forced by the fact that $X$ is smooth) are far stronger.
\proclaim Proposition 3.2.  If the singularity at the generic point of $C$
is a ${{\rm cA}}_2$ singularity, then any dissident singularity $P$ will be
a ${{\rm cD}}_4$ or ${{\rm cE}}_r$ singularity ($r=6,7,8$).  Moreover, the
exceptional divisor $E$ on $X$ has only a pinch point singularity on the fibre
above a dissident point.\par

\noindent \it Proof. \rm  This is just a question of checking the various cases. 
If for instance $P$ were a ${{\rm cA}}_n$ singularity ($n\ge 3$), then we could
choose local analytic coordinates $(x,y,z,t)$ at $P$ so that $C$ is given by
$x=y=z=0$ and  $\Xbar$ is locally given by an equation of the form $x^2 + y^2 +
z^{n+1} + tg(x,y,z,t) =0$.  Our assumptions imply that all terms of $g$ are at
least quadratic in $x,y,z$ and that there is no term in $z^2$ only.  Therefore,
if we blow up the line $x=y=z=0$, we obtain an affine piece of $X$ with equation
$$0=x^2 +y^2 + z^{n-1} + \hbox {\rm other quadratic or higher degree terms in }
x,y,z,t,$$ contradicting the assumption that $X$ was smooth.

We now consider the case of a ${{\rm cD}}_n$ singularity ($n\ge 4$).  Again, we
choose local coordinates so that $C$ is given by
$x=y=z=0$ and  $\Xbar$ is locally given by an equation of the form 
$x^2 + y^2 z + z^{n-1} + tg(x,y,z,t) =0$, where the terms of $g$ are at
least quadratic in $x,y,z$.  If $n>4$, there must be a $z^2$ term in
$g$ (otherwise there is an obvious affine piece of the blow-up which is
singular), but no term in $y^2$ (since for $t\ne 0$ we have a ${{\rm cA}}_2$
singularity).  We then check that a different affine piece of $X$ is
still singular.  In the case when $n=4$, we have the alternative that $g$
contains a term in $y^2$ but no term in $z^2$, and so the equation takes the
form $0=x^2 + y^2 z + z^3 + ty^2 + \hbox {\rm further terms}$. Moreover, by
taking the local equation for $\Xbar$ as  $0=x^2 + y^2 z + z^3 + ty^2$, we see
this case cannot be ruled out by local arguments.  A relevant affine piece of
the blow-up $X$ is then given by an equation $$0 = x^2 + y^2 z + z + ty^2 +
\hbox {terms involving } z,$$ and so $E$  (given locally by $z=0$) has only the
singularity claimed on the dissident fibre.

A similar argument works for the ${{\rm cE}}_r$ singularities ($r=6,7,8$). 
Again we 
choose local coordinates so that $C$ is given by
$x=y=z=0$ and  $\Xbar$ is locally given by an equation of the form 
$h(x,y,z) + tg(x,y,z,t) =0$, where $h = x^2 + y^3 + z^4$ for $r=6$, $h = x^2 +
y^3 + yz^3$ for $r=7$, and  $h = x^2 + y^3 + z^5$ for $r=8$, and 
where the terms of $g$ are at least quadratic in $x,y,z$.  Similar arguments
the those above force the existence of a $z^2$ term in $g$ (and by taking the
local equation for $\Xbar$ to be $0 =h +tz^2$, we see that the case of 
${{\rm cE}}_r$ singularities cannot be ruled out).  Again however, by blowing
up we see that $E$ has only a pinch point singularity (of the form $x^2 + tz^2
= 0$) on the dissident fibre.
\bigskip\bigskip

\centerline{\bf 4. \quad Proof of Main Claim for Type III contractions}\bigskip

\rm Suppose now that we have a codimension one face of the \Kahler cone of
$(X,J)$ corresponding to a Type III contraction $\phimapX$.  Let $E$ denote the
exceptional divisor, which from the results of the previous Section we know to
be a conic bundle over a smooth curve $C$.

We now observe that if $C$ has genus $g>0$, the surface $E$ will not deform
under a generic holomorphic deformation of $(X,J)$, and moreover will only
deform over a locus in the Kuranishi space of codimension at least $g$.  As it
is sufficient to prove it for first order deformations, this can be proved
using a generalization of the Hodge theoretic arguments from [\WKc ].  It will
however also follow from the fact that if $(X,J)$ were to contain such a
conic bundle $E$ over $C$, there would be a Abel-Jacobi map from the Jacobian
$J(C)$ to the intermediate Jacobian $J(X)$, a principally polarized analytic
torus, which is an inclusion (cf.~\S 4 of [\WKc ]).  
Considering the period map, it follows from the 
unobstructedness theorem of Bogomolov, Tian and Todorov that it is a
codimension $g$ condition in moduli for $J(X)$ to contain the $g$-dimensional
abelian variety $J(C)$.

Thus if $g>0$, only finitely many fibres of $E$ will deform under the generic
deformation of $(X,J)$.  We can however also observe that the converse is true,
that if $g=0$ then $E$ will deform for all holomorphic deformations of
$(X,J)$.  This follows from the fact that $E$ is Gorenstein with $\hOE 1 =0$.  
From the exact sequence
$$ 0 \to \calO _X (-E) \to \calO _X \to \calO _E \to 0$$
it follows that $\h i X {\calO _X (E)} =0$ for $i>0$, and hence that 
$\chi (\calO _X (E))=1$.  As argued in \S 2 of [\WKc ], it then follows from the
Base Change Theorem that the divisor class $E$ remains effective when we deform
the complex structure, and so the surface $E$ deforms as claimed.

We concentrate first on the case when the general fibre of $E$ over $C$ is
irreducible (i.e.~isomorphic to $\P ^1$).

When $g>0$, we can apply the arguments from [\WKc ] directly.  If $E$ has a
singular fibre, then the fibre contains a curve $l \cong \P ^1$ with $E\cdot l
=-1$ and $l$ not moving in $X$.  A standard Hilbert scheme argument (see [\WKc ]
\S 3) then shows that $l$ must deform under any holomorphic deformation of the
complex structure on $X$.  If $E$ contains no singular fibres, then it is a
ruled surface over $C$.  A straightforward dimension counting argument
(Proposition 4.2 from [\WKc ]) shows that when $E$ is ruled over a curve of genus
$g>1$, some fibres of $E$ do deform under a holomorphic deformation of $(X,J)$.
Thus, unless $E$ is ruled over an elliptic curve, the given codimension one
face of the \Kahler cone corresponds to a Type I contraction on the generic
holomorphic deformation, contracting finitely many curves isomorphic to $\P
^1$.  The Main Claim therefore follows from the Type I case, as proved in \S 1.

To conclude the argument for the case when $E$ has irreducible generic fibre,
we need to consider the case $g=0$.  This was not a problem in [\WKc ], since as
commented above the whole surface deforms under a holomorphic deformation of
the complex structure on $X$.  This does not however tell us what happens when
we deform to an almost complex structure on $X$.  

The case when $E$ is a ruled surface over $\P ^1$ may be dealt with by the
cobordism argument from [\Ru ] which we outlined at the end of \S 2.  If $A$ is
the class of a fibre of $E$, then the corresponding Hilbert scheme is
isomorphic to $\P ^1$, and so in this case the three point Gromov-Witten
invariant is strictly negative and the Main Claim follows.  We remark in passing
that this argument gives another proof of the Main Claim for the case when $E$
is ruled over a curve $C$ of genus $g>1$, since in that case the Gromov-Witten
invariant  is strictly positive.

For the case when $C$ has genus 0 but $E$ has singular fibres, we need to
replace the Hilbert scheme argument referred to above by a local argument using
the same tools as were used for Type I contractions.  This will of course work
for an arbitrary genus base curve, and so (given the remark in the previous
paragraph) we can if we wish avoid completely the argument of deforming to a
Type~I contraction.  Let $Z$ denote a dissident fibre of $E$ lying over $Q\in
C$, so that $Z$ is either a line pair or a double line.

\proclaim Lemma 4.1.  We can find an open disc $Q\in \Delta \subset C$ and an
open neighbourhood $U$ of $Z$ in $X$ with $U\to \Delta$ a differentiably
trivial family containing the fibre of $E$ over $P$ for all $P\in \Delta$, and
a holomorphic deformation of the complex structure on $U$ for which the
dissident fibre $Z$ deforms into at least two disjoint $(-1,-1)$-curves, but no
other fibre deforms.\par

\noindent \it Proof.  \rm This is very nearly a repetition of the proof of
(1.1).  Let $\phimapX$ be the contraction; then $\Xbar$ is locally just a
family of rational double point surface singularities, parametrized by a disc
$\Delta \subset C$ with centre $Q$.  We therefore have a neighbourhood $U$ of
$Z$ in $X$ with $U \to \Delta$ and fibre $U_0$ over $Q$ containing the dissident
fibre $Z$.

In the notation of [\Fr ], page 678, we have a map $g : \Delta \to S$ so that the
diagram
$$ \matrix{\Delta &\longrightarrow &S\hfill\cr &\searrow \; &\, \downarrow
\pi\hfill\cr && R\hfill\cr }$$
commutes (where $R= \hbox{\rm Def } \bar U_0$ and $S= \hbox{\rm Def } U_0$).

We may now holomorphically deform the family $U \to \Delta$ by means of
deforming $g$, first so that the image of $Q$ is unchanged but no other point
$P\in \Delta$ has image in $D= \pi ^{-1} (\hbox{\rm discriminant
locus})$, and then so that the map is tranverse to $D$.  In this way we obtain 
(shrinking $\Delta$ if necessary) a
deformation $\calU \to \Delta \times \Delta '$ of  $U \to \Delta$ , where for 
$t\ne 0$, the fibre $Z$ deforms to at least
two disjoint $(-1,-1)$-curves, and no other fibre of $E_1 = E|_{\Delta}$ deforms 
(i.e.~we resolve the singularities of the corresponding family $\bar\calU \to
\Delta$, where for $t\ne 0$, the threefold $\bar\calU _t$ has singular locus
consisting of $\delta >1$ nodes). 

Since $\Delta \times \Delta '$ is smooth, we may take a good representative
$\bar \calU$ of the family over $\Delta \times \Delta '$ by Theorem 2.8 of
[\Loo ],
and hence deduce as in (1.1) that $\calU$ is differentiably trivial over 
$\Delta \times \Delta '$.  Taking $U = \calU _0$, the Lemma is proved. \medskip

In order to complete our proof of the Main Claim for the case when the generic
fibre of $E$ over $C$ is irreducible, we need the following corollary of (4.1).

\proclaim Proposition 4.2.  The Main Claim is true in the Type III case
when $E$ has irreducible generic fibre but is not ruled over $C$.\par

\noindent\it Proof.  \rm Let $A$ be the homology class of an irreducible
component of any singular fibre $Z$; we remark that $E\cdot A =-1$ 
and that $A$ is primitive.  As in
(4.1), we can find an open  neighbourhood $U$ of $Z$ with $U \to \Delta$ and
$E_1 =E|_{\Delta}$ having only the one singular fibre $Z$ (above $Q\in \Delta$),
and a holomorphic deformation of the complex structure on $U$ so that $Z$ splits
into at least two $(-1,-1)$-curves.  We may do this for each of the singular
fibres $Z_1 =Z, Z_2, \ldots Z_N$, obtaining deformed holomorphic structures on
each of the corresponding neighbourhoods $U_1 =U,\ldots ,U_N$.  These may be
patched together with the original complex structure $J$ on $X$ to give an
almost complex structure $J'$ on $X$, with respect to which the primitive class
$A$ is represented by a finite number of $(-1,-1)$-curves lying in the open sets 
$U_1,\ldots ,U_N$.  Moreover, as in the proof of (1.1) we may assume that these
are the only pseudo-holomorphic rational curves representing the class $A$.
This latter statement follows as in (1.1) from Gromov compactness,
since by considering the compact subset $F$ of $X$ given by the complement of
suitably small open neighbourhoods of the singular fibres $Z_i$, we may assume
that any other \phrc on $(X,J')$ representing $A$ intersects $F$.  If this were
true for all the almost complex structures $J'$ constructed in the above way
however near to $J$ they were, we would deduce from Gromov compactness that there
was a holomorphic curve on $(X,J)$ representing $A$ and intersecting $F$, which
is not the case.

We remark in passing that it is this step which fails when $E$ is a ruled
surface (say over an elliptic curve $C$).  Here, we have to take $A$ to be the
homology class of a fibre, and then by local holomorphic deformations we can 
always, in a similar way to (4.1),
force certain specified fibres to deform to $(-1,-1)$-curves, which therefore
count positively towards the Gromov-Witten invariant.  We can however no longer
say that these are  the only pseudo-holomorphic rational curves representing the
class $A$, since any fibre of $E$ could be the limit of pseudo-holomorphic
rational curves representing $A$ for $J'$ tending to $J$.  In the case of $C$
elliptic, there must of course be pseudo-holomorphic rational curves on $(X,J')$
which count negatively towards the Gromov-Witten invariant and precisely cancel
out the positive contributions.

In the case under consideration however, we argue as in before that the
Gromov-Witten invariant $\Phi _{A,J'} (H,H,H)$ is strictly positive, and hence 
the Main Claim follows.\bigskip

Having proved the Main Claim under the assumption that the generic fibre of $E$
is irreducible, we assume from now on that the generic fibre of $E$ over $C$ is a
line pair.  In this case, we have a family of lines on $X$ parametrized by a
smooth curve $\tC$, namely the double cover of $C$ branched over the points of
$C$ whose fibres are double lines.  The Hilbert scheme 
corresponding to the lines will be just $\tC$ with embedded components at the
dissident points.  We prove the Main Claim first under the additional
assumption that $g(\tC )>0$.

We observe that a general line in the family will deform over a locus in the
Kuranishi space for $X$ of codimension at most one (a standard argument -- see
(3.1) of [\WKc ]).  Let $\tE$ denote the normalization of $E$ (a ruled surface
over $\tC$) and $f: \tE \to X$ the induced map to $X$.  I claim that the map 
$H^1 (T_X ) \to H^1 (\omega _{\tE})$ (i.e.~$H^{2,1} (X) \to H^{2,1} (\tE )$)
is surjective.  This follows from Hodge theory and the fact that the
contraction $\Xbar$ of $X$ has only quotient singularities in codimension 2, by
an entirely analogous argument to that of (4.1) in [\WKc ].  Observe now that the
above map factors as $H^1 (T_X ) \to H^1 (f^*T_X ) \to H^1
(N_f )\to H^1 (\omega _{\tE})$, and so the obstruction to the morphism $f$
deforming to first order gives at least a codimension $g(\tC )$ condition 
in $H^1(T_X )$, i.e.~$f$ deforms over a locus in the Kuranishi
space of codimension at least $g(\tC )$.  This fact also follows by
consideration of the intermediate Jacobian $J(X)$, since under the above
circumstances there is an Abel-Jacobi inclusion map of the Jacobian $J(\tC )$
into $J(X)$.

So, still under the assumption that $g(\tC )>0$, we have various cases.  If
$g(C)=0$, then $E$ will deform under all complex deformations of $(X,J)$, and
under the generic deformation it must deform to a rational conic bundle with
irreducible generic fibre; this case has already been covered.  Let us assume
therefore that we also have $g(C)>0$.  If $g(\tC )>1$, the
facts stated in the previous paragraph and a dimension counting argument on the
Hilbert scheme of lines in the versal deformation $\calX \to B$ shows that some
line will always deform (cf.~Erratum and (3.1) of [\WKc ]).  Since $E$ does not
deform for the generic deformation, it follows that on such a deformation we
have a Type I contraction, and our Main Claim follows from the results of \S
1.  If however $g(\tC )=1$, then our assumptions imply that $g(C) =1$ and the
map $\tC \to C$ is unramified -- this then is the case of $E$ being a quasi-ruled
surface over an elliptic curve, which because of the results of the Erratum to
[\WKc ] has been excluded when stating the Main Claim.

The only other possibility left is the case that $g(\tC )=0$; here, $C$
is also rational and the double cover is ramified over two points
(i.e.~$E$ has two double fibres).  If $C_0$ denotes the line of
singularities of $E$, the blow-up of $E$ in $C_0$ will just be the
normalization $\tE $, a rational scroll (cf.~(3.2)).  The calculations performed
in the proof of (2.3) apply equally well in this case to show that $C_0$ is a
$(-1,-1)$-curve and that $E$ is isomorphic to the non-normal del Pezzo surface 
$\bar\F _{3,2}$.  In particular, $E$ may be contracted locally analytically to
an isolated Gorenstein singularity.  We now argue as in (2.4) that by suitably deforming the holomorphic structure locally 
on a neighbourhood $U$ of $E$, we may smooth $E$ to a non-singular del Pezzo surface
of degree 7, therefore containing three $(-1,-1)$-curves.  In the Type III case however, 
these curves no longer all have the same numerical class, in that the class corresponding 
to the singular locus of $E$ is different from that of the other two $(-1,-1)$-curves.  If 
we choose our homology class $A$ to be the class of one of these other two 
$(-1,-1)$-curves, then $A$ will define the codimension one face of $\calK$ that 
we want.  Moreover, the argument from (2.4) can be applied with the 
primitive class $A$ to give the 
positivity of the Gromov-Witten invariant $\Phi _A (H,H,H)$, 
and hence to prove the Main Claim in this case.

The Main Claim has therefore now been proved for all the three types of
primitive contraction.
\bigskip\bigskip

\centerline{\bf Acknowledgments} \bigskip\rm
Much of this work was completed by the author whilst a Visiting Scholar for a
month at Harvard University in Spring 1995, and he would like to thank Prof.
S.-T. Yau for making this visit possible.  Thanks are also due to Mark Gross,
Sheldon Katz and Dusa McDuff for valuable conversations and comments at various
stages of the work. 
\bigskip\bigskip

\centerline{\bf References} \bigskip

\rm \noindent [\BW ] \quad J.W. Bruce, C.T.C. Wall, \it On the classification
of cubic surfaces.  \rm J. London Math.~Soc.~19 (1979) 245-256.\par

\noindent [\CFKM ] \quad P. Candelas, A. Font, 
D.R. Morrison, S. Katz, \it Mirror Symmetry for two parameter Models - II, \rm 
preprint hep/9403187, March 1994.\par

\noindent [\Dem ] \quad M. Demazure, \it Surfaces de Del Pezzo. \rm In : M.
Demazure, H. Pinkham,  B. Teissier (eds.), \it
Seminaire sur les singularit\'es des surfaces, \rm pp. 23-70.  Lecture Notes
in Mathematics 777, Springer 1980.\par

\noindent [\Fr ] \quad R. Friedman,  \it  Simultaneous Resolution of
Threefold Double Points, \rm Math. Ann.  274 (1986) 671-689. \par

\noindent [\GrCan ] \quad M. Gross,  \it  The deformation space of \CY
$n$-folds with canonical singularities can be obstructed, \rm Cornell
Univ.~preprint 1994. Duke alg-geom eprint 9402014.\par

\noindent [\GrDef ] \quad M. Gross, \it  Deforming \CY threefolds, \rm
Cornell Univ.~preprint 1995. Duke alg-geom eprint 9506022.\par

\noindent [\HW ] \quad F. Hidaka, K. Watanabe, \it Normal Gorenstein surfaces
with ample anti-canonical divisor, \rm Tokyo J. Math. 4 (1981) 319-330.\par

\noindent [\Kont ] \quad M. Kontsevich, \it Homological Algebra of Mirror
Symmetry, \rm Proc. ICM, Zurich 1994.  Duke alg-geom eprint 9411018.\par

\noindent [\Loo ] \quad E.J.N. Looijenga, \it  Isolated singular points on
complete intersections, \rm  CUP 1984.\par

\noindent [\McD ] \quad D. McDuff,  \it Elliptic methods in symplectic geometry,
 \rm Bull. Amer. Math. Soc. 23 (1990), 311-358.\par

\noindent [\McDS ] \quad D. McDuff, D. Salamon,  \it J-holomorphic curves and
quantum cohomology,  \rm Univ. Lecture Series, Vol. 6,  Amer. Math. Soc.
1994.\par

\noindent [\RCan ]\quad M. Reid, \it  Canonical 3-folds. \rm  In : A. Beauville 
(ed.), \it G\'eom\'etrie alg\'ebrique, Angers 1979, \rm pp. 273-310.  Sijthoff
and Noordhoff, The Netherlands 1990.\par

\noindent [\RMin ]\quad M. Reid, \it  Minimal models of canonical 3-folds. \rm
In : S. Iitaka (ed.), \it Algebraic varieties and analytic varieties, \rm 
Advanced Studies in Pure Math. Vol. 1, pp. 395-418.  Kinokuniya, Tokyo and North
Holland, Amsterdam 1983.\par

\noindent [\ReDP ] \quad M. Reid, \it  Nonnormal del Pezzo surfaces, \rm 
Publications of RIMS, Kyoto, 30:5 (1994) 695-727.\par

\noindent [\Ru ]\quad Y. Ruan, \it Symplectic Topology and extremal rays, \rm
Geometric and Functional Analysis 3 (1993) 395-430. \par 

\noindent [\Ruan ]\quad Y. Ruan, \it  Topological Sigma Model and Donaldson type
invariants in Gromov Theory, \rm to appear in Duke J. Math.\par

\noindent [\RT ]\quad Y. Ruan, G. Tian, \it  A Mathematical Theory of Quantum
Cohomology, \rm to appear in J. Diff. Geom.\par

\noindent [\Sik ]\quad J.-C.~Sikorav, \it Some properties of holomorphic curves
in almost complex manifolds, \rm In : M. Audin, J. Lafontaine (eds), \it
Holomorphic curves in symplectic geometry, \rm Ch.~V, pp.~ 165-190. 
Birkh\"auser Progr. in Math. Vol.~117, Basel-Boston-Berlin 1994.\par

\noindent [\Wall ]\quad C.T.C.~Wall, \it Classification problems in Differential
Topology V.  On certain 6-manifolds, \rm Invent. math. 1 (1966) 355-374.\par

\noindent [\WPic ] \quad P.M.H.~Wilson, \it \CY manifolds with large Picard
number, \rm Invent. math. 98 (1989) 139-155.\par

\noindent [\WKc ] \quad P.M.H.~Wilson, \it The \Kahler cone on \CY
threefolds,  \rm Invent. math. 107 (1992) 561-583. \it Erratum, \rm  
Inventiones math. 114 (1993) 231-233.\par

\noindent [\WEll ] \quad P.M.H.~Wilson, \it Elliptic ruled surfaces on
Calabi-Yau threefolds,  \rm Math. Proc. Cam. Phil. Soc. 112 (1992)
45-52.\par

\bigskip \rightline{DEPARTMENT OF PURE MATHEMATICS}
\rightline {UNIVERSITY OF CAMBRIDGE}\par

\bigskip
\bigskip \it \noindent Mathematics Subject Classification (1991):\hfill\break
14J10, 14J15, 14J30, 32J17, 32J27, 53C15, 53C23, 57R15, 58F05 \hfill\bigskip

\noindent \rm PMHW August 1995, revised February 1996\hfill\break
\noindent e-mail : pmhw@pmms.cam.ac.uk
\hfill\vfill

\hfill\vfill

\end